\documentclass[12pt,a4paper]{article}
\usepackage{amsmath}
\usepackage{amssymb}
\usepackage{amsthm}
\usepackage{amsfonts}
\usepackage{latexsym}

\theoremstyle{plain}
\newtheorem{thm}{Theorem}

\newtheorem{lem}{Lemma}[section]
\newtheorem{prop}{Proposition}[section]
\theoremstyle{definition}

\newtheorem{exmp}{Example}[section]

\theoremstyle{remark}
\newtheorem{claim}{Claim}[section]

\title{Scaling limits for equivariant Szeg\"{o} kernels}
\author{Roberto Paoletti\footnote{\noindent{\bf Address:}
Dipartimento di Matematica e Applicazioni, Universit\`a degli Studi
di Milano Bicocca, Via R. Cozzi 53, 20125 Milano,
Italy; {\bf e-mail}: roberto.paoletti@unimib.it }}
\date{}

\begin{document}

\maketitle

\section{Introduction}

Let $(M,J)$ be an n-dimensional complex projective manifold, and let $(L,h)$ be
an Hermitian ample line bundle on $M$. Suppose that the unique compatible
connection on $L$ has curvature $\Theta=-2 i\,\omega$, where $\omega$ is a
Hodge form on $M$. The pair $(\omega,J)$ puts an Hermitian structure
$H=g-i\omega$ on the (holomorphic) tangent bundle
$TM$, hence a Riemannian structure $g$ on $M$.

Let $G$ be a compact
connected g-dimensional Lie group, and suppose given a Hamiltonian holomorphic
action of $G$ on $(M,\omega,J)$ unitarily linearizing to $(L,h)$.
For every $k=1,2,\ldots$, there is a natural Hermitian structure
on each space of holomorphic global sections
$H^0(M,L^{\otimes k})$, and a naturally induced
unitary representation of $G$ on $H^0(M,L^{\otimes k})$.

Let $\{V_\varpi\}_{\varpi \in \Theta}$ be the finite dimensional irreducible
representations of $G$, and for every $\varpi\in \Theta$ let $H^0(M,L^{\otimes k})_\varpi
\subseteq H^0\left (M,L^{\otimes k}\right )$ be the maximal subspace equivariantly isomorphic
to a direct sum of copies of $V_\varpi$. There are unitary equivariant isomorphisms
\begin{equation}\label{eqn:action-on-sections}
H^0(M,L^{\otimes k})=\bigoplus _{\varpi \in \Theta}H^0(M,L^{\otimes k})_\varpi.\end{equation}
The action of $G$ on $L$ dualizes to an action on the dual line bundle $L^*$
in a natural manner; on the other hand,
the $G$-invariant Hermitian metric $h$ on $L$ naturally induces an Hermitian metric
on $L^*$, still denoted by $h$, which is also $G$-invariant.

Let $X\subseteq L^*$ is the unit circle bundle, with projection
$\pi:X\rightarrow M$. Then by the above the action of $G$ on $L^*$ leaves $X$ invariant.
Furthermore, $X$ is a contact manifold, with contact form given by the connection
1-form $\alpha$. Since $G$ preserves both the Hermitian metric and the holomorphic
structure, it preserves the unique compatible connection, and therefore it acts on $X$
as a group of contactomorphisms; given this, $X$ has a standard
$G$-invariant Riemannian metric. By these underlying structures,
in the following we shall tacitly identify functions, densities and
half-densities on $X$. In the following, to avoid cumbersome notation, we shall use the same symbol $\mu_g$
for the symplectomorphism of $M$ and the contactomorphism of $X$ induced by $g\in G$.

As is well-known, the spaces of smooth sections $\mathcal{C}^\infty\left (M,L^{\otimes k}\right )$
may be unitarily and equivariantly identified with the spaces $\mathcal{C}^\infty(X)_k$
of smooth functions on $X$ of the $k$-th isotype for the $S^1$-action, that is, obeying the covariance law
$f(e^{i\vartheta}\cdot x)=e^{ik\vartheta}\,f(x)$ for $x\in X$ and $e^{i\vartheta}\in S^1$.
Let $H(X)_k\subseteq \mathcal{C}^\infty(X)_k$ be the subspace of functions corresponding to $H^0(M,L^{\otimes k})$
under this isomorphism, so that $H(X)=:\bigoplus _{k=0 }^{+\infty}H(X)_k$ is the Hardy space of $X$.
Thus (\ref{eqn:action-on-sections}) translates into

\begin{equation}\label{eqn:action-on-functions}
H(X)_k=\bigoplus _{\varpi \in \Theta}H(X)_{\varpi,k}.\end{equation}

In this paper, we are concerned with certain $\mathcal{C}^\infty$ functions $\Pi _{\varpi,k}$
on $X$ naturally
associated to each pair $(\varpi,k)\in \Theta\times \mathbb{N}$. Namely, let us choose for
any $(\varpi,k)\in \Theta\times \mathbb{N}$ an orthonormal basis
$\left \{s^{(\varpi,k)}_j\right \}_{j=1}^{N_{\varpi k}}$ of $H(X)_{\varpi,k}$, and let us define
$$
\Pi_{\varpi,k}(x,y)=:\sum _{j=1}^{N_{\varpi k}}s^{(\varpi,k)}_j(x)\cdot \overline{s^{(\varpi,k)}_j(y)}
\,\,\,\,\,\,\,\,\,(x,y\in X).
$$
Then $
\Pi_{\varpi,k}$ is well-defined, that is, independent of the choice of the orthonormal basis, and in fact it can be
intrinsically described as the distributional kernel of the orthogonal projection
$P _{\varpi,k}:L^2(X)\rightarrow H(X)_{\varpi,k}$. We shall study here the asymptotic properties of the functions
$
\Pi_{\varpi,k}$, as $\varpi$ is fixed and $k\rightarrow +\infty$.

Let $\mathfrak{g}$ be the Lie algebra of $G$, and denote by $\Phi:M\rightarrow \mathfrak{g}^*$
the moment map of the action of $G$ on $(M,2\omega)$. In \cite{pao-mm},
it has been shown that for fixed
$\varpi$ one has
$\Pi_{\varpi,k}(x,x)=O(k^{-\infty})$ as $k\rightarrow +\infty$, unless $\Phi \big (\pi(x)\big)=0$.

On the other hand, if
$\Phi \big (\pi(x)\big)=0$,
and $G$ acts freely on $\Phi ^{-1}(0)\subseteq M$, then by Corollary 1 of \cite{pao-sq}
(working with a different
normalization convention for the total volume) there is an asymptotic expansion
$$
\Pi_{\varpi,k}(x,x)=\sum _j\left |s^{(\varpi,k)}_j(x)\right |^2\sim \frac{\dim (V_\varpi)^2}{V_{\mathrm{eff}}\big(x\big)}\,
k^{\mathrm{n}-\mathrm{g}/2}+\sum _{l\ge 1}a_{l,\varpi}(x)\,k^{\mathrm{n}-\mathrm{g}/2-l},$$
where $V_{\mathrm{eff}}:(\Phi\circ \pi)^{-1}(0)\rightarrow \mathbb{R}$ is the effective potential of the action \cite{burnsg};
its value
on $x\in (\Phi\circ \pi)^{-1}(0)$ is the
volume of the $G$-orbit in $M$ through $\pi(x)$.
Thus the effective potential of the action controls the asymptotics of the restriction of $
\Pi_{\varpi,k}$ to the diagonal of $X\times X$.

In the particular case of the trivial representation $\varpi =0$,
$V_{\mathrm{eff}}$ relates the asymptotics of $\Pi_{0,k}$ and of the Szeg\"{o} kernel of the symplectic reduction
$(M_0,\omega_0,L_0)$ of $(M,L,\omega)$, expressing an obstruction to the conformal
unitarity of the Guillemin-Sternberg
map $H^0\left (M,L^{\otimes k}\right )^G\rightarrow H^0\left(M_0,L_0^{\otimes k}\right )$. Further developments on this
problem are due to Charles \cite{ch}, Hall and Kirwin \cite{hk}, Hui Li \cite{l}, Ma and Zhang \cite{mz}.

Turning momentarily to the action free case, the fast decay of Szeg\"{o} kernels
away from the diagonal has stimulated interest in the asymptotics of their scaling limits near the diagonal.
More precisely, suppose $x\in X$, and let $\rho(z,\theta)$
be a Heisenberg local chart for $X$ centered at $x$,
as in (\ref{eqn:heisenberg-coordinates}) below; in particular, if $m=:\pi(x)$
this unitarily
identifies $(T_mM,H_m)$ and $\mathbb{C}^\mathrm{n}$ with its standard Hermitian structure.
As shown in Theorem 3.1 of \cite{sz},
for any $w,v\in \mathbb{C}^\mathrm{n}$
the following asymptotic expansion holds as $k\rightarrow +\infty$
for the level-$k$ Szeg\"{o} kernel $\Pi_k$:
\begin{eqnarray}\label{eqn:expansion-for-pi}
%
\Pi_{k}\left
(\rho\left(\frac{v}{\sqrt{k}},\theta\right),\rho\left(\frac{w}{\sqrt{k}},\theta'\right)\right)
\sim\left
(\frac{k}{\pi}\right )^{\mathrm{n}}\, e^{ik(\theta-\theta')+\psi _2(u,v)}\,\left (1+\sum _{j\ge
1}a_j(x,w,v)\,k^{-j/2}\right),
\end{eqnarray}
where
$$\psi _2(w,v)=:w\cdot \overline{v}-\frac 12\,\left (\|w\|^2+\|v\|^2\right
),
$$ and the $a_j$ are polynomials in $w$ and $v$ (see
also \cite{bsz} for the leading term).
Recall that here $H_m$ denotes the Hermitian structure of $TM$ induced by
$\omega =:\frac i2 \,\Theta$; this normalization convention accounts for the
factor $\frac{1}{\pi^\mathrm{n}}$ in
(\ref{eqn:expansion-for-pi}), unlike the earlier work \cite{z}.
We shall conform here to \cite{sz}; thus the total volume of
$M$ is $\mathrm{vol}(M)=\frac{\pi^\mathrm{n}}{\mathrm{n}!}\,\int _Mc_1(L)^\mathrm{n}$.

In this article, we shall study the scaling limits of the
equivariant Szeg\"{o} kernels $\Pi _{\varpi,k}$, and show that to
leading order they are still simply related to the effective volume,
certain data associated to the representation $\varpi$, and (in the
special case where $G$ acts freely on $\Phi^{-1}(0)\subseteq M$) the
scaling limits of the Szeg\"{o} kernel of the symplectic reduction.
Furthermore, we shall see that equivariant scaling limits can also be expressed by
the product of an exponentially decaying factor in $v,\,w$ times an
asymptotic expansion whose coefficients are polynomials in $v$ and $w$.
We remark that in the toric case equivariant asymptotics have been studied in \cite{stz}.

To express our results, we need some basic facts about the local
geometry of $M$ along $M'=:\Phi ^{-1}(0)$ \cite{gs-gq}, \cite{ggk}.
Recall that if $0\in \mathfrak{g}^*$ is a regular value of $\Phi$,
then $M'$ is a g-codimensional connected coisotropic submanifold of
$M$, whose null-fibration is given by the orbits of the $G$-action.

At any $m\in M$, let us denote by $\mathfrak{g}_M(m)\subseteq T_mM$
the tangent space to the orbit through $m$, and by
$J_m:T_mM\rightarrow T_mM$ the complex structure.

If $m\in M'$, let us denote by $Q_m\subseteq T_mM$ the Riemannian
orthocomplement of $\mathfrak{g}_M(m)$ in $T_mM'$. Thus, $Q_m$ is a
complex subspace of $T_mM'$, of complex dimension
$\mathrm{n}-\mathrm{g}$.

The Riemannian orthocomplement of $T_mM'\subseteq T_mM$ is $J_m\Big
(\mathfrak{g}_M(m)\Big)$. Therefore, we have orthogonal direct sum
decompositions

\begin{equation}\label{eqn:direct-sum}
T_mM=T_mM'\oplus J_m\Big (\mathfrak{g}_M(m)\Big),\,\,\,\,
T_mM'=Q_m\oplus \mathfrak{g}_M(m).
\end{equation}
Given (\ref{eqn:direct-sum}), if $m\in M'$ and $w\in T_mM$, we shall
decompose $w$ as $w=w_{\mathrm{v}}+w_{\mathrm{h}}+w_{\mathrm{t}}$,
where $w_{\mathrm{v}}\in \mathfrak{g}_M(m)$, $w_{\mathrm{h}}\in
Q_m$, $w_{\mathrm{t}}\in J_m\Big (\mathfrak{g}_M(m)\Big )$. The
labels stand for vertical, horizontal, and transverse. This hints to
the fact that in the special case where $G$ acts freely on $M'$, the
latter is a principal $G$-bundle on the symplectic reduction
$M_0=M'/G$; thus $\mathfrak{g}_M(m)$ is the vertical tangent fibre,
while $Q$ is a connection projecting unitarily to the tangent bundle
of $M_0$.

Before stating our Theorem, another definition is in order. To this
end, recall that if $0\in \mathfrak{g}^*$ is a regular value of the
moment map then the action of $G$ on $\Phi ^{-1}(0)\subseteq M$ is
locally free. Therefore, any $m\in \Phi ^{-1}(0)$ has finite
stabilizer subgroup $G_m\subseteq G$.

Suppose $x\in X$, $\Phi \big (\pi(x)\big )=0$.
If $G_{\pi(x)}\subseteq G$ is the (finite) stabilizer subgroup of $\pi(x)$,
for every $g\in G_{\pi(x)}$ there exists a unique $h_g\in S^1$ such that
$\mu _g(x)=h_g\cdot x$, where $\mu _g:X\rightarrow X$ is the contactomorphism induced by $g$. We shall then let
\begin{equation}\label{eqn:key-function}
A_{\varpi,k}(x)=:2^{\mathrm{g}/2}\,\frac{\dim (V_\varpi)}{V_{\mathrm{eff}}(x)}\,\cdot
\frac{1}{\left |G_{\pi(x)}\right |}\,\sum _{g\in G_{\pi(x)}}\chi _\varpi(g)\,h_g^{k},
\end{equation}
where $\chi _\varpi:G\rightarrow \mathbb{C}$ is the character of the
irreducible representation $\varpi$.

As above, $\omega =:\frac i2 \,\Theta$, where $\Theta$ is the curvature, and $h$ is the Hermitian metric on $TM$ associated
to $\omega$.

Furthermore, as in (\ref{eqn:expansion-for-pi}) we shall express the asymptotic expansion for
$\Pi_{\varpi,k}$ in a Heisenberg local chart $\rho$ centered at $x$.
However, given that the dependence of $\Pi_{\varpi,k}$ on $\theta$ and $\theta'$ is given by the factor
$e^{ik(\theta-\theta')}$ and carries no geometric information, in the
following we shall generally take $\theta=\theta'=0$; with the identification
$T_mM\cong \mathbb{C}^\mathrm{d}$ induced by $\rho$ understood, we shall set $x+w/\sqrt{k}=:\rho\big(w/\sqrt{k},0\big)$.

We then have:

\begin{thm} \label{thm:equiv-scaling-limit}
Suppose that $0\in \mathfrak{g}^*$ is a regular value of $\Phi$, and
$x\in X$, $\Phi\big (\pi(x)\big)=0$. Let us choose a system of Heisenberg local coordinates
centered at $x$. For every $\varpi \in \Theta$ and $w,v\in
T_{\pi(x)}M$, the following asymptotic expansion holds as
$k\rightarrow +\infty$:
\begin{eqnarray*}
\lefteqn{\Pi_{\varpi,k}\left
(x+\frac{w}{\sqrt{k}},x+\frac{v}{\sqrt{k}}\right)}\\
 &&\sim \left
(\frac{k}{\pi}\right
)^{\mathrm{n}-\mathrm{g}/2}\,A_{\varpi,k}(x)\,e^{Q(w_\mathrm{v}+w_\mathrm{t},
v_{\mathrm{v}}+v_\mathrm{t})}\,e^{\psi
_2(w_\mathrm{h},v_\mathrm{h})}\cdot \left (1+\sum _{j\ge 1}a_{\varpi
j}(x,w,v)\,k^{-j/2}\right),\end{eqnarray*} where
$$
Q(w_\mathrm{v}+w_\mathrm{t},
v_{\mathrm{v}}+v_\mathrm{t})=-\|v_\mathrm{t}\|^2-\|w_\mathrm{t}\|^2+i\,\big[\omega
_m(w_\mathrm{v},w_\mathrm{t})-\omega
_m(v_\mathrm{v},v_\mathrm{t})\big],$$
and the $a_{\varpi
j}$'s are polynomials in $v$, $w$ with coefficients depending on $x$ and $\varpi$.

\end{thm}

We integrate the statement by the following remarks.

\begin{itemize}
%
  \item The remainder term can be given a \lq large ball estimate\rq $\,$(that is, for
$\|u\|,\,\|v\|\lesssim k^{1/6}$), similar to the ones in \cite{sz}.
More precisely, let $R_N(x,v,w)$ be the remainder term following the first $N$ summands in
(\ref{eqn:asympt-expansion-general-case}).
Given the description of
$\Pi_{\varpi,k}$ as an oscillatory integral (cfr (\ref{eqn:j-th-summand-fourier-itegral-exponential-in-t-theta}
below), we may adapt the arguments in \S 5 of \cite{sz}
to obtain that for $\|u\|,\|v\|\lesssim k^{1/6}$ we have
$$
\left |R_N(x,v,w)\right |\le C_N\,k^{\mathrm{n}-(\mathrm{g}+N+1)/2}\,e^{
-\frac{1-\epsilon}{2}\,
\left(\|u_\mathrm{h}-v_\mathrm{h}\|^2+2\|v_\mathrm{t}\|^2+2\|w_\mathrm{t}\|^2\right)}.
$$
The bound also holds in $\mathcal{C}^j$-norm.
  \item In the special case where $G$ acts freely on
$\Phi^{-1}(0)$, denote by $X_0\subseteq L_0^*$ the circle bundle
of the reduced pair $(M_0,L_0)$, and by $\Pi_{k}^{(0)}$ the level
$k$ Szeg\"{o} kernel of $X_0$. If $\Phi \big (\pi(x)\big )=0$, let
us denote by $\overline{x}$ its image in $X_0$, and if
$w_\mathrm{h}\in Q_{\pi(x)}$ let $\overline{w}_\mathrm{h}$ be its
isometric image in the tangent space to $M_0$. By
(\ref{eqn:expansion-for-pi}) and Theorem
\ref{thm:equiv-scaling-limit}, we obtain

\begin{eqnarray*}
\lefteqn{\Pi_{\varpi,k}\left
(x+\frac{w}{\sqrt{k}},x+\frac{v}{\sqrt{k}}\right)}\\
 &\sim&2^{\mathrm{g}/2}\,\left
(\frac{k}{\pi}\right )^{\mathrm{n}-\mathrm{g}/2}\,\frac{\dim
(V_\varpi)^2}{V_{\mathrm{eff}}(x)}\,e^{Q(w_\mathrm{v}+w_\mathrm{t},
v_{\mathrm{v}}+v_\mathrm{t})}\,e^{\psi
_2(w_\mathrm{h},v_\mathrm{h})}+\sum _{j\ge 1}a_{\varpi
j}(x,w,v)\,k^{\mathrm{n}-(\mathrm{g}+j)/2}\\
&=&\left(\frac{2k}{\pi}\right)^{\mathrm{g}/2}\cdot \left(\frac{\dim
(V_\varpi)^2}{V_{\mathrm{eff}}(x)}\,e^{Q(w_\mathrm{v}+w_\mathrm{t},
v_{\mathrm{v}}+v_\mathrm{t})}\right )\cdot \Pi_{k}^{(0)}\left
(\overline{x}+\frac{\overline{w}_\mathrm{h}}{\sqrt{k}},\overline{x}+\frac{\overline{v}_\mathrm{h}}{\sqrt{k}}\right)
+\mathrm{L.O.T.}.\end{eqnarray*}
  \item Arguing as in \S 2.3 of \cite{dp}, one can see that
$\Pi_{\varpi,k}=O\left(k^{-\infty}\right)$ uniformly on compact subsets of the complement
in $X\times X$ of the locus
$$I(\Phi)=:\left\{(x,y):x\in (G\times S^1)\cdot y,\,\Phi\big(\pi(y)\big)=0\right\}.$$
Thus it is natural to consider scaling limits at any $(x,y)\in I(\Phi)$.
Given $g_0\in G$, $h_0\in S^1$, $x\in \left (\Phi\circ \pi\right)^{-1}(0)$ and $v,w\in T_{\pi(x)}(M)$,
a minor modification of the arguments in the proof of Theorem \ref{thm:equiv-scaling-limit} leads to
an asymptotic expansion
\begin{eqnarray}
\label{eqn:asympt-expansion-general-case}
\lefteqn{\Pi_{\varpi,k}\left
(\mu_{g_0}\circ r_{h_0}\left(x+\frac{w}{\sqrt{k}}\right),x+\frac{v}{\sqrt{k}}\right)}\\
 &&\sim \left
(\frac{k}{\pi}\right
)^{\mathrm{n}-\mathrm{g}/2}\,A_{\varpi,k}(x,g_0,h_0)\,e^{Q(w_\mathrm{v}+w_\mathrm{t},
v_{\mathrm{v}}+v_\mathrm{t})}\,e^{\psi
_2(w_\mathrm{h},v_\mathrm{h})}\cdot \left (1+\sum _{j\ge 1}a_{\varpi
j}(x,w,v)\,k^{-j/2}\right),\nonumber\end{eqnarray}
where now
$$
A_{\varpi,k}(x,g_0,h_0)=:2^{\mathrm{g}/2}\,\frac{\dim (V_\varpi)}{V_{\mathrm{eff}}(x)}\,\cdot
\frac{1}{\left |G_{\pi(x)}\right |}\,\sum _{g\in G_{\pi(x)}}\chi _\varpi \left(g\,g_0^{-1}\right)\cdot
(h_0\,h_g)
^{k}.
$$
\item We are primarily interested in the case of complex projective
manifolds. In view of the microlocal description of almost complex
Szeg\"{o} kernels appearing in \cite{sz}, the results of this paper
can however be extended to the context of almost complex symplectic
manifolds.

\end{itemize}

After this paper was completed, I learned of the rich paper
\cite{mz} alluded to above.
Using analytic
localization techniques of Bismut and Lebeau for $\mathrm{spin}^c$ Dirac operators, Ma and Zhang
obtain among other things an asymptotic
expansion for the trivial representation.

\noindent
\textbf{Acknowledgments.} I am indebted to Steve Zelditch for a
remark that improved the statement of Theorem \ref{thm:equiv-scaling-limit}, and to the referee for suggesting various improvements
in presentation.

\section{Examples}


In the non-equivariant case, a key feature of scaling asymptotics of Szeg\"{o} kernels
expressed by (\ref{eqn:expansion-for-pi})
is the universal nature of the leading term,
essentially the level-one Szeg\"{o} kernel of the reduced Heisenberg group $\mathbf{H}^\mathrm{n}_{\mathrm{red}}$.
To express this more precisely, recall that the latter may be viewed as the unit circle bundle
of the trivial line bundle $L=\mathbb{C}^\mathrm{n}\times \mathbb{C}$ over $\mathbb{C}^\mathrm{n}$, endowed with the Hermitian
metric
$$
h\big((\mathbf{z},w),(\mathbf{z},w')\big)=w\,\overline{w'}\,e^{-\|\mathbf{z}\|^2}\,\,\,\,\,\,\Big(
\mathbf{z},\mathbf{z}'\in \mathbb{C}^\mathrm{n},\,w,w'\in \mathbb{C}\Big).
$$
The unit circle bundle is thus given by
$$
X=\mathbf{H}^\mathrm{n}_{\mathrm{red}}=\left\{(\mathbf{z},w)\in \mathbb{C}^\mathrm{n}\times \mathbb{C}:
|w|=e^{-\|\mathbf{z}\|^2/2}\right\}.
$$
A Heisenberg chart for $X$ centered at $(\mathbf{0},1)$ is
$$
\varphi_\mathbf{0}:\mathbb{C}^\mathrm{n}\times (-\pi,\pi)\rightarrow X,\,\,\,\,\,\,\left(\mathbf{z},\theta\right)
\mapsto \left (\mathbf{z},e^{-\|\mathbf{z}\|^2/2+i\theta}\right).
$$
As shown in \cite{bsz}, for every $k=1,2,\ldots$ the level-k Szeg\"{o} kernel is
\begin{equation}
\label{eqn:level-k-heisenberg}
\Pi_k^\mathbf{H}\big(\varphi_\mathbf{0}(\mathbf{w},\theta),\varphi_\mathbf{0}(\mathbf{v},\theta')\big)=
\left(\frac{k}{\pi}\right)^\mathrm{n}\,e^{k\big[i(\theta-\theta')+\psi_2(\mathbf{w},\mathbf{v})\big]}.
\end{equation}

In the linear case, we shall derive
from (\ref{eqn:level-k-heisenberg})
an asymptotic expansion in the spirit of Theorem \ref{thm:equiv-scaling-limit}, at any $x=\left(\mathbf{z_1},e^{-\mathbf{z_1}/2}\right)$
for which the map $\gamma_\mathbf{z_1}:g\in G\mapsto \mu_g(\mathbf{z_1})\in \mathbb{C}^\mathrm{d}$ is an embedding
(that is, $\mathbf{z}_1$ has trivial stabilizer in $G$); with minor changes,
the arguments below apply when $\gamma_\mathbf{z_1}$ is an immersion (that is,
$\mathbf{z}_1$ has finite stabilizer in $G$).


\begin{exmp}
Let $A:G\rightarrow \mathrm{U}(\mathrm{n}),\,g\mapsto A_g$, be a unitary representation, so that the underlying
action on $\left (\mathbb{C}^\mathrm{n},\omega_0\right)$ is $\mu_g(\mathbf{z})=:A_g\mathbf{z}$ ($\mathbf{z}\in
\mathbb{C}^\mathrm{n}$); here $\omega_0=:\frac{i}{2}\,\sum_{j=1}^\mathrm{n}dz_j\wedge d\overline{z}_j$
is the standard symplectic structure on $\mathbb{C}^\mathrm{n}$.
The standard Hermitian structure on $\mathbb{C}^n$ is then
$H_0=g_0-i\omega_0$, where $g_0(\mathbf{w},\mathbf{v})=-\omega_0(J_0\mathbf{w},\mathbf{v})$
($J_0$ being multiplication by $i$).

A linearization to $L$
is given by
$$
\mu_g\big((\mathbf{z},w)\big)=:\big(A_g\mathbf{z},w\big).
$$

For any $\mathbf{z}_1\in \mathbb{C}^\mathrm{n}$,
a Heisenberg chart for $X$ centered at
$\left (\mathbf{z}_1,e^{-\|\mathbf{z}_1\|^2/2}\right)$ is
$$
\varphi_{\mathbf{z}_1}:\left (\mathbf{z},\theta\right)\mapsto \varphi_\mathbf{0}
\big(\mathbf{z}+\mathbf{z}_1,\omega_0(\mathbf{z},\mathbf{z}_1)+\theta\big)=
\left(\mathbf{z}+\mathbf{z}_1,
e^{-\|\mathbf{z}+\mathbf{z}_1\|^2/2+i\big(\omega_0(\mathbf{z},\mathbf{z}_1)+\theta\big)}
\right).
$$
Thus, given $x=\left (\mathbf{z}_1,e^{-\|\mathbf{z}_1\|^2/2}\right)\in \mathbf{H}^\mathrm{n}_{\mathrm{red}}$
and $\mathbf{v}\in \mathbb{C}^\mathrm{n}$, in our notation
$$
x+\mathbf{v}=\varphi_{\mathbf{z}_1}(\mathbf{v},0)=\varphi_\mathbf{0}\big(\mathbf{z}_1+\mathbf{v},
\omega_0(\mathbf{v},\mathbf{z}_1)\big).
$$

Given an irreducible representation $\varpi$ and $\mathbf{w},\mathbf{v}\in \mathbb{C}^\mathrm{n}$, by a straightforward
computation using (\ref{eqn:level-k-heisenberg}) we obtain
\begin{eqnarray}
\label{eqn:espressione-integrale-per-pi}
\lefteqn{\Pi^{\mathbf{H}}_{\varpi,k}\left (x+\frac{\mathbf{w}}{\sqrt{k}},x+\frac{\mathbf{v}}{\sqrt{k}}\right)=\dim(V_\varpi)\cdot}
\nonumber\\
&&\int_G\chi_\varpi(g)\,\Pi^\mathbf{H}_k\left(
\varphi_\mathbf{0}\left(
A_g\mathbf{z}_1+\frac{A_g\mathbf{w}}{\sqrt{k}},\frac{1}{\sqrt{k}}\,\omega_0\left(\mathbf{w},\mathbf{z}_1\right)
\right),
\varphi_\mathbf{0}\left(
\mathbf{z}_1+\frac{\mathbf{v}}{\sqrt{k}},\frac{1}{\sqrt{k}}\,\omega_0\left(\mathbf{v},\mathbf{z}_1\right)\right)\right)\,dg
\nonumber\\
&=&\dim(V_\varpi)\,\left(\frac{k}{\pi}\right)^\mathrm{n}\,\int_G\chi_\varpi(g)\,
e^{S_k(\mathbf{z}_1,\mathbf{w},\mathbf{v})}\,dg,
\end{eqnarray}
where $dg$ is the density on $G$ associated to an invariant Riemannian metric of total volume one, and
\begin{eqnarray}
\label{eqn:fase-for-S-k}
S_k(\mathbf{z}_1,\mathbf{w},\mathbf{v})
&=:&
k\,H_0\big(A_g\mathbf{z}_1-\mathbf{z}_1,\mathbf{z}_1\big)\\
&&+\sqrt{k}\,\Big[H_0\left (\mathbf{v},A_g^{-1}\mathbf{z}_1-\mathbf{z}_1\right)
+H_0\left (A_g\mathbf{z}_1-\mathbf{z}_1,\mathbf{w}\right)\Big]+
\psi_2(A_g\mathbf{v},\mathbf{w}).\nonumber
\end{eqnarray}

Given the simplifying assumption that $\mathbf{z}_1$ has trivial stabilizer in $G$,
there exists $C>0$ such that $\|A_g\mathbf{z}_1-\mathbf{z}_1\|\ge C\,\mathrm{dist}_G(g,e)$,
where $e\in G$ is the unit and $\mathrm{dist}_G$ is the Riemannian metric on $G$.
Thus, it follows from (\ref{eqn:level-k-heisenberg}) that the integrand of (\ref{eqn:espressione-integrale-per-pi})
is $O\left(k^{-\infty}\right)$ on the loci $A_k\subseteq G$ where, say,
$\mathrm{dist}_G(g,e)\ge 2\,k^{-1/3}$.
On the loci $B_k$ where $\mathrm{dist}_G(g,e)\le k^{-1/3}$, on the other hand, we can transfer the integration
to the Lie algebra $\mathfrak{g}$ by the exponential map $\exp_G:\mathfrak{g}\rightarrow G$,
$\eta\mapsto e^\eta$, and apply the rescaling $\eta=\frac{1}{\sqrt{k}}\,\xi$.
Let $\widehat{A}:\mathfrak{g}\rightarrow \mathfrak{u}(\mathrm{n})$, $\eta\mapsto \widehat{A}_\eta$,
be the differential of the
morphism of Lie groups $A:G\rightarrow \mathrm{U}(\mathrm{n})$.
Thus
$$
A_{e^{\xi/\sqrt{k}}}=e^{\widehat{A}_{\xi/\sqrt{k}}}=\mathrm{id}_{V}+\frac{\widehat{A}(\xi)}{\sqrt{k}}+\frac{1}{2}\,\frac{\widehat{A}(\xi)^2}{k}+
\cdots.
$$
On the upshot, after some computations we obtain
\begin{eqnarray}
\label{eqn:integration-over-g-linear}
\lefteqn{\Pi_{\varpi,k}\left(x+\frac{w}{\sqrt{k}},x+\frac{v}{\sqrt{k}}\right)
\sim \frac{\dim(V_\varpi)}{\pi^\mathrm{n}}\,k^{\mathrm{n}-\mathrm{g}/2}\,\int_{\mathfrak{g}}
\chi_\varpi\left (e^{\xi/\sqrt{k}}\right)\,e^{\mathcal{S}_k(\xi,\mathbf{w},\mathbf{v},\mathbf{z}_1)}\,d\xi}\nonumber\\
&=&\frac{\dim(V_\varpi)^2}{\pi^\mathrm{n}}\,k^{\mathrm{n}-\mathrm{g}/2}\,\int_{\mathfrak{g}}
e^{\mathcal{S}_k(\xi,\mathbf{w},\mathbf{v},\mathbf{z}_1)}\,d\xi\cdot \left (1+O\left(k^{-1/2}\right)\right),
\end{eqnarray}
where now
\begin{eqnarray}
\label{eqn:phase-for-linear-case}
\lefteqn{\mathcal{S}_k(\xi,\mathbf{w},\mathbf{v},\mathbf{z}_1)}\nonumber\\
&=:&
i\sqrt{k}\,\omega_0\left(\mathbf{z}_1,\widehat{A}_\xi(\mathbf{z}_1)\right)+\psi_2(\mathbf{w},\mathbf{v})-
\frac 12\,\left \|\widehat{A}_\xi(\mathbf{z}_1)\right\|^2+H_0\left (\widehat{A}_\xi(\mathbf{z}_1),\mathbf{v}\right)
-H_0\left (\mathbf{w},\widehat{A}_\xi(\mathbf{z}_1)\right)\nonumber\\
&=&
i\sqrt{k}\,\Phi^\xi\big(\mathbf{z}_1\big)+\psi_2(\mathbf{w},\mathbf{v})-
\frac 12\,\left \|\widehat{A}_\xi(\mathbf{z}_1)\right\|^2+H_0\left (\widehat{A}_\xi(\mathbf{z}_1),\mathbf{v}\right)
-H_0\left (\mathbf{w},\widehat{A}_\xi(\mathbf{z}_1)\right);
\end{eqnarray}
here $\Phi:V\rightarrow \mathfrak{g}^*$ is the moment map, and $\Phi^\xi=:\langle \Phi,\xi\rangle$.

Suppose, to begin with, that $\Phi\big(\mathbf{z}_1\big)\neq \mathbf{0}$. Then the linear phase
$\xi\mapsto \Phi^\xi\big(\mathbf{z}_1\big)$ has no stationary point in $\xi$, and since by (\ref{eqn:phase-for-linear-case})
the integrand in (\ref{eqn:integration-over-g-linear}) is absolutely convergent, the stationary phase Lemma applies
to show that $\Pi_{\varpi,k}\left(x+\frac{w}{\sqrt{k}},x+\frac{v}{\sqrt{k}}\right)=
O\left(k^{-\infty}\right)$.

If $\Phi\big(\mathbf{z}_1\big)= \mathbf{0}$, on the other hand, we have
\begin{eqnarray*}
\lefteqn{\Pi_{\varpi,k}\left(x+\frac{w}{\sqrt{k}},x+\frac{v}{\sqrt{k}}\right)}\\
&=&\frac{\dim(V_\varpi)^2}{\pi^\mathrm{n}}\,k^{\mathrm{n}-\mathrm{g}/2}\,
e^{\psi_2(\mathbf{w},\mathbf{v})}\,\int _{\mathfrak{g}}e^{-
\frac 12\,\left \|\widehat{A}_\xi(\mathbf{z}_1)\right\|^2+H_0\left (\widehat{A}_\xi(\mathbf{z}_1),\mathbf{v}\right)
-H_0\left (\mathbf{w},\widehat{A}_\xi(\mathbf{z}_1)\right)}\,d\xi\cdot \left (1+k^{-1/2}\right)\\
&=&\frac{\dim(V_\varpi)^2}{\pi^\mathrm{n}\cdot V_{\mathrm{eff}}(\mathbf{z}_1)}\,k^{\mathrm{n}-\mathrm{g}/2}\,
e^{\psi_2(\mathbf{w},\mathbf{v})}\,\int _{\mathfrak{g}(\mathbf{z}_1)}e^{-
\frac 12\,\left \|\mathbf{s}\right\|^2+H_0\left (\mathbf{s},\mathbf{v}\right)
-H_0\left (\mathbf{w},\mathbf{s}\right)}\,d\mathbf{s}\cdot \left (1+k^{-1/2}\right),
\end{eqnarray*}
where in the latter equality integration has been shifted from $\mathfrak{g}$ to the tangent space
$\mathfrak{g}(\mathbf{z}_1)\subseteq \mathbb{C}^\mathrm{n}$ at $\mathbf{z}_1$ to the $G$-orbit of $\mathbf{z}_1$ by the change of variables
$\mathbf{s}= \widehat{A}_\xi(\mathbf{z}_1)$; hence $d\mathbf{s}=V_{\mathrm{eff}}(\mathbf{z}_1)\,d\xi$.

Given the equalities
$H_0\big(\mathbf{w},\mathbf{s}\big)=H_0\big(\mathbf{w}_\mathrm{t}+\mathbf{w}_\mathrm{v},\mathbf{s}\big)$,
$\psi_2\big(\mathbf{w},\mathbf{v}\big)=\psi_2\big(\mathbf{w}_{\mathrm{h}},\mathbf{v}_{\mathrm{h}}\big)
+\psi_2\big(\mathbf{w}_{\mathrm{t}}+\mathbf{w}_{\mathrm{v}},\mathbf{v}_{\mathrm{t}}+\mathbf{v}_{\mathrm{v}}\big),
$
we obtain with a few calculations
\begin{eqnarray}\label{eqn:gaussian-integral-linear}
\lefteqn{\Pi_{\varpi,k}\left(x+\frac{w}{\sqrt{k}},x+\frac{v}{\sqrt{k}}\right)}\\
&=&\frac{\dim(V_\varpi)^2}{\pi^\mathrm{n}\cdot V_{\mathrm{eff}}(\mathbf{z}_1)}\,k^{\mathrm{n}-\mathrm{g}/2}\,
e^{\psi_2(\mathbf{w},\mathbf{v})+\frac 12\|\mathbf{w}_\mathrm{v}-\mathbf{v}_\mathrm{v}\|^2}\,
\int _{\mathfrak{g}(\mathbf{z}_1)}e^{-i\omega_0(\mathbf{s},
\mathbf{v}_\mathrm{t}+\mathbf{w}_\mathrm{t})-
\frac 12\,\left \|\mathbf{s}-(\mathbf{w}_\mathrm{v}-\mathbf{v}_\mathrm{v})\right\|^2}\,d\mathbf{s}\cdot \left (1+O\left(k^{-1/2}\right)\right).
\nonumber
\end{eqnarray}
The Gaussian integral in (\ref{eqn:gaussian-integral-linear}) is
$
(2\pi)^{\mathrm{g}/2}\,e^{i\omega_0(\mathbf{v}_\mathrm{t}+\mathbf{w}_\mathrm{t},\mathbf{w}_\mathrm{v}-\mathbf{v}_\mathrm{v})-\frac 12\|\mathbf{v}_\mathrm{t}
+\mathbf{w}_\mathrm{t}\|^2}$, and from this one computes
\begin{eqnarray*}
\lefteqn{\Pi_{\varpi,k}\left(x+\frac{w}{\sqrt{k}},x+\frac{v}{\sqrt{k}}\right)}\\
&=&2^{\mathrm{g}/2}\,\frac{\dim(V_\varpi)^2}{ V_{\mathrm{eff}}(\mathbf{z}_1)}\,\left(\frac k\pi\right)^{\mathrm{n}-\mathrm{g}/2}\,
e^{\psi_2(\mathbf{w}_\mathrm{h},\mathbf{v}_\mathrm{h})-\|\mathbf{w}_\mathrm{t}\|^2-\|\mathbf{v}_\mathrm{t}\|^2
+i\big(\omega_0(\mathbf{w}_\mathrm{v},\mathbf{w}_\mathrm{t})-\omega_0(\mathbf{v}_\mathrm{v},\mathbf{v}_\mathrm{t})\big)}
\cdot \left (1+O\left(k^{-1/2}\right)\right).
\nonumber
\end{eqnarray*}
\end{exmp}

Before considering the next example, let us recall from \cite{bsz} that for $k=1,2,\ldots$ an orthonormal basis of
$H^0\left(\mathbb{P}^\mathrm{d},\mathcal{O}_{\mathbb{P}^\mathrm{d}}(k)\right)$
is $\left\{s_\mathbf{J}^k\right\}_{|\mathbf{J}|=k}$, where
\begin{equation}
\label{eqn:orthonormal-basis-general}
s_\mathbf{J}^k=:\sqrt{\frac{(k+\mathrm{d})!}{\pi^\mathrm{d}\,\mathbf{J}!}}\,z^\mathbf{J};
\end{equation}
here $\mathbf{J}!=:\prod _{l=0}^\mathrm{d}j_l!$,
$z^\mathbf{J}=:\prod _{l=0}^\mathrm{d}z_l^{j_l}$.

\begin{exmp} The unitary representation of $S^1$ on $\mathbb{C}^2$ given by
$t\cdot (z_0,z_1)=:\left (t^{-1}z_0,tz_1\right)$ descends to a symplectic action
on $\mathbb{P}^1$, with a built-in linearization to the hyperplane line bundle.
The associated moment map is
$$\Phi:\mathbb{P}^1\rightarrow \mathbb{R},\,\,\,\,\,\,[z_0:z_1]\mapsto \frac{-|z_0|^2+|z_1|^2}{|z_0|^2+|z_1|^2}.$$

Clearly, any $[z_0:z_1]\in \Phi^{-1}(0)$ has stabilizer subgroup $\{\pm 1\}$. Since any $S^1$-orbit
in $S^3$ has length $2\pi$ and doubly covers its image in $\mathbb{P}^1$, the effective volume
is identically equal to $\pi=\frac{2\pi}{2}$ on $\Phi^{-1}(0)$.
Therefore,
$$
A_{\varpi,k}\big([z_0:z_1]\big)=\frac{\sqrt{2}}{\pi}\cdot\frac 12\,\left[1+(-1)^\varpi\,(-1)^k\right]=
\left\{\begin{array}{ccc}
         \sqrt{2}/\pi & \mathrm{if} & k\equiv \varpi\,(\mathrm{mod}.\,2) \\
         0 & \mathrm{if} & k\not\equiv \varpi\,(\mathrm{mod}.\,2).
       \end{array}\right.
$$

Given that
$$
\mu_t\left (z_0^l\,z_1^{k-l}\right)=\left(z_0^l\,z_1^{k-l}\right)\circ \mu_{t^{-1}}= \left(t z_0\right)^l\,\left(t^{-1}z_1\right)^{k-l}=t^{2l-k}\,z_0^l\,z_1^{k-l},
$$
we have for $\varpi\in \mathbb{Z}$ and $k\in \mathbb{N}$:
\begin{equation}
\label{eqn:orthonormal-basis-case-of-p1}
H^0\left(\mathbb{P}^1,\mathcal{O}_{\mathbb{P}^1}(k)\right)_\varpi=\left\{
\begin{array}{ccc}
  \mathrm{span}\left\{z_0^{\frac{\varpi +k}{2}}\,z_1^{\frac{k-\varpi}{2}}\right\} & \mathrm{if} & k\equiv \varpi\,\,(\mathrm{mod.}\,2), \\
 0& \mathrm{if} &k\not\equiv \varpi\,\,(\mathrm{mod.}\,2).
\end{array}\right.
\end{equation}
By the Stirling formula, if $b$ is fixed and $a\rightarrow +\infty$ we have
\begin{eqnarray}
\label{eqn:stirling}
(a+b)!
\sim\sqrt{2\pi\,a}\,\left(\frac{a^{a+b}}{e^a}\right).
\end{eqnarray}
Suppose then $k=\varpi+2s$, $s\in \mathbb{N}$, and choose $(z_0,z_1)\in S^3$ lying over $[z_0:z_1]$; in view of (\ref{eqn:orthonormal-basis-general}), (\ref{eqn:orthonormal-basis-case-of-p1})
and (\ref{eqn:stirling}),
\begin{eqnarray}
\label{eqn:stirling-applied}
\Pi_{\varpi,k}\big([z_0:z_1],[z_0:z_1]\big)&=&
\frac{(\varpi+2s+1)!}{\pi\,(\varpi+s)!\,s!}\,|z_0|^{2(\varpi+s)}\,|z_1|^{2s}\\
&\sim&\frac 1\pi\,\sqrt{\frac s\pi}\,2^{\varpi+2s+1}\,|z_0|^{2(\varpi+s)}\,|z_1|^{2s},\nonumber
\end{eqnarray}
as $s\rightarrow +\infty$.

If $[z_0:z_1]\in \Phi^{-1}(0)$, so that $|z_0|^2=|z_1|^2=\frac 12$,
we obtain
\begin{equation*}
\Pi_{\varpi,k}\big([z_0:z_1],[z_0:z_1]\big)\sim \frac 2\pi\,\sqrt{\frac s\pi}\sim
\frac {\sqrt{2}}{\pi}\,\sqrt{\frac{k}{\pi}}=A_{\varpi,k}\big([z_0:z_1]\big)\,\sqrt{\frac{k}{\pi}},
\end{equation*}
which
fits with the asymptotic expansion of Theorem \ref{thm:equiv-scaling-limit}.

If $|z_0|\neq |z_1|$, (\ref{eqn:stirling}) is rapidly decreasing as $s\rightarrow +\infty$.

\end{exmp}

\section{Preliminaries}
\label{sctn:prelimiaries}
In this section we shall collect some preliminaries and set
some notation.

If $(M,J)$ is a complex manifold, any K\"{a}hler form $\omega$ on it
determines an Hermitian metric $h$ on the tangent bundle of $M$, and
$\omega=-\Im(h)$. The Riemannian metric $g=:\Re(h)$ is
$g_m(w,v)=\omega _m(w,J_mv)$ ($m\in M$, $w,v\in T_mM$).

Since Heisenberg local coordinates centered at a given $x\in X$ will
be a key tool in the following, we shall briefly recall their
definition \cite{sz}.

Thus we now assume that $L\rightarrow M$ is an Hermitian ample line
bundle, and $\omega=\frac i2\,\Theta$, where $\Theta$ is the curvature of the
unique compatible covariant derivative.
Let us choose an \textit{adapted} holomorphic
coordinate system $(z_1,\cdots,z_\mathrm{n})$ for $M$ centered at
$\pi(x)$. This means that, when expressed in the $z_i$'s, $\omega$
evaluated at $\pi(x)$ is the standard symplectic structure on
$\mathbb{C}^\mathrm{n}$, that is, $\omega \big (\pi(x)\big)=\frac
i2\,\sum _{j=1}^\mathrm{n}dz_j\wedge d\overline{z}_j$. Thus the
choice of the $z_i$'s determines a unitary isomorphism
$T_{\pi(x)}M\cong \mathbb{C}^\mathrm{n}$.

Let us next choose a \textit{preferred} local frame $e_L$ for $L$ at
$\pi(x)$, in the sense of \cite{sz}. Thus $e_L$ is a holomorphic
local section for $L$ in the neighborhood of $\pi(x)$, satisfying
\begin{equation}\label{eqn:condizioni}
\left \|e_L\big(\pi(x)\big)\right \|=1,\,\,\,\nabla e_L\big
(\pi(x)\big)=0,\,\,\,\nabla ^2e_L\big
(\pi(x)\big)=-\overline{h}_{\pi(x)}\otimes e_L\big (\pi(x)\big),
\end{equation}
where $\nabla$ is the covariant derivative of the connection, and
$\overline{h}=g+i\omega$. The local holomorphic frame for $L$
uniquely determines a holomorphic dual local frame $e_L^*$ for
$L^*$, determined by the condition $\left (e_L^*,e_L\right)=1$,

For $\delta>0$, let $B_{2\mathrm{n}}(0;\delta)\subseteq \mathbb{C}^{\mathrm{n}}
\cong\mathbb{R}^{2\mathrm{n}}$ be the ball of radius $\delta$
centered at the origin.
For an appropriate $\delta >0$, a system of Heisenberg local
coordinates for $X$ centered at $x$ is then given by the map
\begin{equation}
\label{eqn:heisenberg-coordinates}
\rho: B_{2\mathrm{n}}(0;\delta)\times (-\pi,\pi)\rightarrow
X,\,\,\,\, (z,\theta)\mapsto e^{i\theta}\,a(z)^{-1/2}\,e_L^*(z),
\end{equation}
where $a(z)=:\|e_L^*\|^2=\|e_L\|^{-2}$.
If $w\in T_{\pi(x)}M\cong \mathbb{C}^\mathrm{n}$, we shall denote by
$x+w$ the point in $X$ with Heisenberg local coordinates $(w,0)$.

It will simplify our exposition to make a little equivariant adjustment to
the previous construction. Suppose that $m\in M$ has finite stabilizer subgroup
$G_m\subseteq G$ (this will always be the case when $\Phi(m)=0$ if $0\in \mathfrak{g}^*$
is a regular value of the moment map).
Let $U\subseteq M$ be a $G_m$-invariant open neighborhhod of the identity, and suppose
that a local holomorphic frame $\sigma=e_L^*$ satisfying (\ref{eqn:condizioni}) has been chosen on $U$.
Clearly, for every $g\in G_m$ we have $g^*(\sigma)(m)=h_g\cdot \sigma(m)$ (recall that
$g^*(\sigma)=\mu_g\circ\sigma\circ\mu_{g^{-1}}$).
We may then consider the new frame
$$
\overline{\sigma}=\frac{1}{|G_m|}\,\sum _{g\in G_m}\,h_g^{-1}\,g^*(\sigma).$$
Then $\overline{\sigma}(m)=e_L(m)$, and since the metric and the connection are $G$-invariant
$\overline{\sigma}$ also satisfies (\ref{eqn:condizioni}).
Moreover, we now have
\begin{equation}\label{eqn:equiv-preferred}g^*(\overline{\sigma})=h_g\cdot \overline{\sigma},
\,\,\,\,\forall\,g\in G_m.\end{equation}
In the following, the underlying preferred local holomorphic frame
in the definition of Heisenberg local coordinates will be assumed to satisfy (\ref{eqn:equiv-preferred}).

For $\xi \in \mathfrak{g}$, we shall denote by $\xi _M$ and $\xi _X$
the vector fields on $M$ and $X$, respectively, associated to $\xi$
by the infinitesimal actions of $\mathfrak{g}$. The moment map
$\Phi:M\rightarrow \mathfrak{g}^*$ for the action on $(M,2\omega)$
is related to the $G$-invariant
connection form $\alpha$ on $X$ by the relation $\Phi ^\xi =-\iota
(\xi_X)\,\alpha$, where $\Phi^\xi=\left <\Phi,\xi\right >$.

\section{Proof of Theorem \ref{thm:equiv-scaling-limit}.}

To begin with, let us fix an invariant Haar metric on $G$, and let
$dg$ denote the associated measure; by Haar metric we mean that
$\int _Gdg=1$. Now if $\rho:G\rightarrow \mathrm{GL}(W)$ is
linear representation on a complex vector space, for any $\varpi \in
\Theta$ the projection $P_\varpi$ of $W$ onto the the
$\varpi$-isotypical component $W_\varpi$ is given by
\begin{equation}
\label{eqn:projection-varpi} P_\varpi =\dim (V_\varpi)\,\int _G\chi
_\varpi\left (g^{-1}\right )\,\rho (g)\,dg
\end{equation}
\cite{dixmier}. On the other hand, the unitary representation of $G$ on
$H_k(X)\subseteq L^2(X)$ induced by the action on $X$ is given by
$(g\cdot f)(y)=:f\left (\mu_{g^{-1}}(y)\right )$ ($f\in L^2(X)$, $y\in
X$). Therefore, the equivariant Szeg\"{o} kernel $\Pi _{\varpi,k}$
is given by
\begin{equation}\label{eqn:equivariant-varpi-k-szego-integral}
\Pi_{\varpi,k}(y,y')=\dim (V_\varpi)\,\int _G\chi
_\varpi\left(g^{-1}\right )\,\Pi_k\left (\mu_{g^{-1}}(y),y'\right
)\,dg,
\end{equation}
where $\mu_g:X\rightarrow X$ is the contactomorphism associated to
$g\in G$.

Suppose $x\in X$, $\Phi (x)=0$, and set $m=:\pi(x)$. We assume given
a system of Heisenberg local coordinates for $X$ centered at $x$.
This choice gives a meaning to the expression $x+w$, for any $w\in
T_mM\cong \mathbb{C}^\mathrm{n}$.

Then for every $\varpi\in \Theta$ and $k\in \mathbb{N}$ we have
\begin{eqnarray}\label{eqn:proj-da-k-omega}
\lefteqn{\Pi_{\varpi,k}\left
(x+\frac{w}{\sqrt{k}},x+\frac{v}{\sqrt{k}}\right)}\\
&&=\dim(V_\varpi)\,\int_G\,\chi_\varpi
\left(g^{-1}\right)\,\Pi_k\left(\mu_{g^{-1}}\left(x+\frac{w}{\sqrt{k}}\right),
x+\frac{v}{\sqrt{k}}\right)\,dg,\nonumber\end{eqnarray}
where $\chi _\varpi:G\rightarrow \mathbb{C}$ is the character of the
irreducible representation $V_\varpi$ \cite{dixmier}.

We shall now split the integration in $d\mu $ as the sum of two
terms, one which is rapidly decaying as $k\rightarrow +\infty$, and
another where integration is over a suitably shrinking neighborhood
of the (finite) stabilizer subgroup $G_m\subseteq G$.

To this end, let us define for every $k\in \mathbb{N}$ an open cover $\{A_k,\, B_k\}$ of $G$ by
setting
$$A_k=:\left \{g\in G:\mathrm{dist}_G\big (g,G_m\big)>k^{-1/3}\right \},$$
$$B_k=:\left \{g\in G:\mathrm{dist}_G\big (g,G_m\big)<2\,k^{-1/3}\right \}.$$
(towards application of the stationary phase Lemma
later in the proof, the exponent $-1/3$ used in the definition of $A_k$ and $B_k$,
could be replaced by $-a$, for any $a\in (0,1/2)$).
Here $\mathrm{dist}_G:G\times G\rightarrow \mathbb{R}$ is the Riemannian distance function. Let $a_k+b_k=1$
be a $G_m$-invariant partition of unity on $G$ subordinate
to the open cover $\{A_k,\, B_k\}$
(thus, $\mathrm{supp}(a_k)\subseteq A_k$, $\mathrm{supp}(b_k)\subseteq B_k$).

We may then split (\ref{eqn:proj-da-k-omega}) as
\begin{eqnarray}\label{eqn:proj-da-k-omega-split}
\lefteqn{\Pi_{\varpi,k}\left
(x+\frac{w}{\sqrt{k}},x+\frac{v}{\sqrt{k}}\right)}\\
&&=\Pi_{\varpi,k}\left
(x+\frac{w}{\sqrt{k}},x+\frac{v}{\sqrt{k}}\right)_a+\Pi_{\varpi,k}\left
(x+\frac{w}{\sqrt{k}},x+\frac{v}{\sqrt{k}}\right)_b;\nonumber\end{eqnarray}
the first (respectively, second) summand in (\ref{eqn:proj-da-k-omega-split})
is (\ref{eqn:proj-da-k-omega}) with the integrand multiplied by $a_k$ (respectively, $b_k$).

\begin{prop}\label{prop:a-term-decays}
$\Pi_{\varpi,k}\left
(x+\frac{w}{\sqrt{k}},x+\frac{v}{\sqrt{k}}\right)_a=O(k^{-\infty})$ as $k\rightarrow +\infty$.
\end{prop}

\textit{Proof of Proposition \ref{prop:a-term-decays}.}
Let $\mathrm{dist}_M:M\times M\rightarrow M$ be the Riemannian distance function.
We have:

\begin{lem}
\label{lem:a-term-decays}
There exists a positive constant $C$
(dependent on $w$ and $v$, but independent of $k$) such that for all $k\gg 0$ and $g\in A_k$
we have
\begin{equation}
\label{eqn:distance-not-small}
\mathrm{dist}_M\left (\mu_{g^{-1}}\left (m+\frac{w}{\sqrt{k}}\right ),m+\frac{v}{\sqrt{k}}\right)
\ge C\,\mathrm{dist}_G(g,G_m)\ge C\,k^{-1/3}.
\end{equation}
\end{lem}

\textit{Proof of Lemma \ref{lem:a-term-decays}.}
If not, we can find $\mathbb{N}\ni k_j\uparrow +\infty$ and
$g_j\in A_{k_j}$ such that $\forall j=1,2,\ldots$ we have
\begin{equation}
\label{eqn:bad-sequence}
\mathrm{dist}_M\left (\mu_{g_j}\left (m+\frac{w}{\sqrt{k_j}}\right ),m+\frac{v}{\sqrt{k_j}}\right)
\le \frac 1j\,\mathrm{dist}_G(g_j,G_m).
\end{equation}
Since $\mathrm{dist}_G(g_j,G_m)$ is bounded above by the diameter of the compact Lie group $G$,
we have in particular
$\mathrm{dist}_M\left (\mu_{g_j}\left (m+\frac{w}{\sqrt{k_j}}\right ),m+\frac{v}{\sqrt{k_j}}\right)
\rightarrow 0$, hence also
$\mathrm{dist}_M\left (\mu_{g_j}\left (m\right ),m\right)\rightarrow 0$.
Therefore, $g_j\rightarrow G_m$; after passing to a subsequence, therefore, we may assume that
$g_j\rightarrow g_0$ for some $g_0\in G_m$. Let us write $g_j= g_0\,h_j$,
where $h_j\rightarrow e$, and $\mathrm{dist}_G(h_j,e)=\mathrm{dist}_G(g_j,G_m)\ge k_j^{-1/3}$.
Using the exponential map $\exp_G:\mathfrak{g}\rightarrow G$, for all $j\gg 0$
we can write
$h_j=e^{\nu_j}$, for unique $\nu_j\in \mathfrak{g}$ such that
$\|\nu_j\|=\mathrm{dist}_G(h_j,e)$.
Since $G$ acts locally freely on $\Phi^{-1}(0)$, there exists $c>0$ such that
$
\|\nu_M(m)\|\ge c\,\|\nu\|,\,\forall\,m\in \Phi^{-1}(0),\,
\nu\in \mathfrak{g}$ (the former norm is in $T_mM$, the latter in
$\mathfrak{g}$). Hence,
\begin{equation}
\label{eqn:bound-on-nu-j-m}
\big\|\big(\nu_j\big)_M(m)\big\|\ge c\,k_j^{-1/3}\,\,\,\,\,\,(j\gg 0).
\end{equation}

Working in preferred local coordinates centered at $m$, we have
\begin{eqnarray}
\label{eqn:action-in-heis-loc-coord}
\mu_{e^{\nu_j}}\left (m+\frac{w}{\sqrt{k_j}}\right )&=&m+\big(\nu_j\big)_M(m)+\frac{w}{\sqrt{k_j}}+O\left(k_j^{-2/3}\right),
\nonumber\\
\mu_{g_0^{-1}}\left(m+\frac{v}{\sqrt{k_j}}\right)&=&m+O\left(k_j^{-1/2}\right).
\end{eqnarray}
By definition of preferred local coordinates, it follows from (\ref{eqn:bound-on-nu-j-m}) and (\ref{eqn:action-in-heis-loc-coord})
that for $j\gg 0$ we have
$$
\mathrm{dist}_M\left (\mu_{e^{\nu_j}}\left (m+\frac{w}{\sqrt{k_j}}\right ),\mu_{g_0^{-1}}\left(m+\frac{v}{\sqrt{k_j}}\right)\right)
\ge \frac c2\, \|\nu_j\|.$$
On the other hand,
we can rewrite (\ref{eqn:bad-sequence}) as
\begin{equation}
\label{eqn:bad-sequence-is}
\mathrm{dist}_M\left (\mu_{e^{\nu_j}}\left (m+\frac{w}{\sqrt{k_j}}\right ),\mu_{g_0^{-1}}\left(m+\frac{v}{\sqrt{k_j}}\right)\right)\le \frac 1j\,\|\nu_j\|,\end{equation}
a contradiction.

\hfill Q.E.D.

\medskip

Returning to the proof of Proposition \ref{prop:a-term-decays},
by Lemma \ref{lem:a-term-decays}
and the off-diagonal estimate on the Szeg\"{o} kernel in (6.1) of \cite{christ}, we conclude
\begin{equation}
\label{eqn:rapid-decay}
\left | \Pi_k\left (\mu_{g^{-1}}\left (m+\frac{w}{\sqrt{k}}\right ),m+\frac{v}{\sqrt{k}}\right)\right |
\le C\,e^{-C_2\,k^{1/6}}
\end{equation}
whenever $k\gg 0$ and $g\in A_k$. The statement follows easily from (\ref{eqn:rapid-decay}).

\hfill Q.E.D.

\medskip

Since our focus is on asymptotic expansions, we shall henceforth disregard the $a$ term.
Let us set $\beta _{\varpi}(g)=\frac{\dim (V_\varpi)}{2\pi}\,\chi_\varpi\left (g^{-1}\right)$ ($g\in G$).
Then

We have
\begin{eqnarray}
\label{eqn:integral-for-b-term}
\lefteqn{\Pi_{\varpi,k}\left
(x+\frac{w}{\sqrt{k}},x+\frac{v}{\sqrt{k}}\right)_b}\\
&=&\int _{-\pi}^\pi\,\int _{B_k}\,\beta _{\varpi}(g)\,\,b_k(g)\,e^{-ik\vartheta}\,
\Pi\left (\mu _{g^{-1}}\circ r_{e^{i\vartheta}}
\left (x+\frac{w}{\sqrt{k}}\right ),x+\frac{v}{\sqrt{k}}\right )\,d\vartheta\,dg\nonumber.
\end{eqnarray}

Suppose $G_m=\{g_1=e,g_2,\ldots,g_{N_x}\}$. Let us define
$$E_k=:\big\{g\in G:\mathrm{dist}_G(g,e)<2\,k^{-1/3}\big\}.$$
Then $B_k=\bigcup _{j=1}^{N_x}B_{jk}$, where $B_{jk}=g_j\cdot E_k$
($1\le j\le N_x$, $k\in \mathbb{N}$). Thus

\begin{eqnarray}
\label{eqn:integral-for-b-term-split-j}
\lefteqn{\Pi_{\varpi,k}\left
(x+\frac{w}{\sqrt{k}},x+\frac{v}{\sqrt{k}}\right)_b}\\
&=&\sum _j\int _{-\pi}^\pi\,\int _{B_{jk}}\,\beta _{\varpi}(g)\,\,b_k(g)\,e^{-ik\vartheta}\,
\Pi\left (\mu _{g^{-1}}\circ r_{e^{i\vartheta}}
\left (x+\frac{w}{\sqrt{k}}\right ),x+\frac{v}{\sqrt{k}}\right )\,d\vartheta\,dg\nonumber\\
&=&\sum _j\Pi_{\varpi,k}\left
(x+\frac{w}{\sqrt{k}},x+\frac{v}{\sqrt{k}}\right)_{j},\nonumber
\end{eqnarray}
where the $j$-th summand in (\ref{eqn:integral-for-b-term-split-j}) is
\begin{eqnarray}\label{eqn:j-th-summand}
\lefteqn{\Pi_{\varpi,k}\left
(x+\frac{w}{\sqrt{k}},x+\frac{v}{\sqrt{k}}\right)_{j}}\\
&=&e^{-ik\vartheta_j}\int _{-\pi}^{\pi}\,\int _{E_k}\,\beta _{\varpi}(g_j\,g)\,\,b_k(g_jg)\,e^{-ik\vartheta}\,
\Pi\left (\mu _{g^{-1}g_j^{-1}}\circ r_{e^{i(\vartheta+\vartheta_j)}}
\left (x+\frac{w}{\sqrt{k}}\right ),x+\frac{v}{\sqrt{k}}\right )\,d\vartheta\,dg;\nonumber
\end{eqnarray}
here $e^{i\vartheta_j}=h_{g_j}$ for every $j$.
Notice that $b_k(g_jg)=b_k(g)$ for every $k\in \mathbb{N}$ and $j$, since $b_k$ is $G_m$-invariant.

Let us examine the asymptotics of (\ref{eqn:j-th-summand}).
To this end, fix $\epsilon>0$ very small (but independent of $k$), and
let $\gamma _0+\gamma _1=1$ be a partition of unity on $(-\pi,\pi)$ with $\mathrm{supp}(\gamma _0)\subseteq
(-\epsilon,\epsilon)$,  $\mathrm{supp}(\gamma _1)\subseteq
(-\pi,-\epsilon/2)\cup (\epsilon/2,\pi)$.
Then
$$
\Pi_{\varpi,k}\left
(x+\frac{w}{\sqrt{k}},x+\frac{v}{\sqrt{k}}\right)_j=\sum _{l=0}^1\Pi_{\varpi,k}\left
(x+\frac{w}{\sqrt{k}},x+\frac{v}{\sqrt{k}}\right)_{jl},
$$
where $\Pi_{\varpi,k}\left
(x+\frac{w}{\sqrt{k}},x+\frac{v}{\sqrt{k}}\right)_{jl}$ is given by (\ref{eqn:j-th-summand})
with the integrand multiplied by $\gamma _l$.

\begin{lem}\label{lem:b1-term-decay-fast}
$\Pi_{\varpi,k}\left
(x+\frac{w}{\sqrt{k}},x+\frac{v}{\sqrt{k}}\right)_{j1}=O(k^{-\infty})$ as $k\rightarrow +\infty$.
\end{lem}

\textit{Proof.}
If $k\gg 0$, $g\in B_{jk}$ and $|\vartheta|>\epsilon /2$, then
$$
\mathrm{dist}_X\left(\mu _{g^{-1}g_j^{-1}}\circ r_{e^{i(\vartheta+\vartheta_j)}}
\left (x+\frac{w}{\sqrt{k}}\right ),x+\frac{v}{\sqrt{k}}\right )>\frac{\epsilon}{3}
$$
($w$ and $v$ are held fixed). Since the singular support of the Szeg\"{o} kernel $\Pi$ is the
diagonal $\mathrm{diag}(X)\subseteq X\times X$, we conclude that
$$
\Psi _{k,g}(h)=:\gamma _0(h)\,\beta _\varpi(g)\,b_k(g)\,\Pi\left(\mu _{g^{-1}g_j^{-1}}\circ r_{e^{i(\vartheta+\vartheta_j)}}
\left (x+\frac{w}{\sqrt{k}}\right ),x+\frac{v}{\sqrt{k}}\right )$$
is a bounded family of smooth functions on $S^1$ when
$k\ge k_0$, $g\in B_{jk}$ and $|\vartheta|>\epsilon /2$; here $\gamma _0$ is interpreted as $\gamma _0 (e^{i\vartheta})$,
a cut-off function supported
on a small open neighborhood of $1\in S^1$.

In the same range, therefore, for every $l\in \mathbb{N}$ we can find a constant
$C_l>0$ such that $\left|\Psi _{k,g}^{(s)}\right |<C_l\,s^{-l}$ for every $s\in \mathbb{N}$, where
$\Psi _{k,g}^{(s)}$ denotes the $s$-th Fourier coefficient of $\Psi _{k,g}$. In particular, this is true for
$s=k$, hence $\left|\Psi _{k,g}^{(k)}\right |<C_l\,k^{-l}$. The same estimate then holds
after integrating over $B_k$, and this implies the statement.

\medskip

We are reduced to studying the asymptotics of $\Pi_{\varpi,k}\left
(x+\frac{w}{\sqrt{k}},x+\frac{v}{\sqrt{k}}\right)_{j1}$. To proceed, let us
introduce the parametrix for the Szeg\"{o} kernel contructed in \cite{bs}. Thus, up to
a smoothing term which does not contribute to the asymptotic expansion, we can
represent $\Pi$ as a Fourier integral operator of the form
\begin{equation}\label{eqn:fourier-integral-for-pi}
\Pi(y,y')=\int _0^{+\infty} e^{it\psi (y,y')}\,s(y,y',t)\,dt\,\,\,\,\,\,(y,y'\in X),
\end{equation}
where the phase satisfies $\Im (\psi)\ge 0$, and the amplitude
is a semiclassical symbol admitting an asymptotic expansion $s(y,y',t)\sim
\sum _{t=0}^{+\infty} t^{\mathrm{n}-j}\,s_j(y,y')$. In view of Lemma \ref{lem:b1-term-decay-fast},
inserting (\ref{eqn:fourier-integral-for-pi})
into (\ref{eqn:j-th-summand}), and multiplying the integrand by $\gamma _0$, we obtain

\begin{eqnarray}\label{eqn:j-th-summand-fourier-itegral}
\lefteqn{\Pi_{\varpi,k}\left
(x+\frac{w}{\sqrt{k}},x+\frac{v}{\sqrt{k}}\right)_{j}}\\
&\sim&e^{-ik\vartheta_j}\int_0^{+\infty} \int _{-\epsilon}^{\epsilon}\,\int _{E_k}\,\gamma_0(\vartheta)\,\beta _{\varpi}(g_j\,g)\,b_k(g)\,e^{i\left [t\psi \left (\mu _{(g_jg)^{-1}}\circ r_{e^{i(\vartheta+\vartheta_j)}}
\left (x+w/\sqrt{k}\right ),x+v/\sqrt{k}\right )-k\vartheta\right] }\nonumber \\
&&\cdot s\left (\mu _{g^{-1}g_j^{-1}}\circ r_{e^{i(\vartheta+\vartheta_j)}}
\left (x+\frac{w}{\sqrt{k}}\right ),x+\frac{v}{\sqrt{k}},t\right )\,dt\,d\vartheta\,dg\nonumber\\
&=&k\,e^{-ik\vartheta_j}\int_0^{+\infty} \,\int _{-\epsilon}^{\epsilon}\,\int _{E_k}e^{ik \Psi _{kj}
(g,t,\vartheta)}\,A_{\varpi kj}(g,t,\vartheta)\,dt\,d\vartheta\,dg;\nonumber
\end{eqnarray}
in the last equality we have performed the coordinate change $t\rightsquigarrow kt$, and set
\begin{equation}\label{eqn:psi-j-k}
\Psi _{kj}
(g,t,\vartheta)=:t\psi \left (\mu _{(g_jg)^{-1}}\circ r_{e^{i(\vartheta+\vartheta_j)}}
\left (x+w/\sqrt{k}\right ),x+v/\sqrt{k}\right )-\vartheta,
\end{equation}

\begin{eqnarray}\label{eqn:a-j-k}
\lefteqn{A_{\varpi kj}(g,t,\vartheta)}\\
&=:&\gamma_0(\vartheta)\,\beta _{\varpi}(g_j\,g)\,b_k(g)\cdot
s\left (\mu _{g^{-1}g_j^{-1}}\circ r_{e^{i(\vartheta+\vartheta_j)}}
\left (x+\frac{w}{\sqrt{k}}\right ),x+\frac{v}{\sqrt{k}},kt\right )\nonumber.
\end{eqnarray}

Let $\exp_G:\mathfrak{g}\rightarrow G$ be the exponential map, and let $E\subseteq \mathfrak{g}$ be
a suitably small open neighborhood of the origin $0\in \mathfrak{g}$, over which $\exp_G$ restricts to a diffeomorphism
$E\rightarrow E'=:\exp_G(E)$. Since
the shrinking open neighborhood $E_k\subseteq G$ of the identity $e\in G$ is definitely contained in
$E'$, we may express the integration in $dg$ using the exponential chart.
To this end, let us fix an orthonormal basis of $\mathfrak{g}$, so as to unitarily identify
$\mathfrak{g}$ with $\mathbb{R}^\mathrm{g}$, and let us write $\xi$ for the correspondig linear coordinates
on $\mathfrak{g}$. We shall denote by $H_G(\xi)\,d\xi$ the local coordinate expression of the
Haar measure $dg$ under the exponential cart; the orthonormality of the chosen basis of $\mathfrak{g}$
implies that $H_G(0)=1$.

With some abuse of language, we shall write $b_k$ for the composition $b_k\circ \exp_G$, and assume that
$b_k(\xi)=b\left (\sqrt[3]{k}\,\xi\right )$ for a fixed function $b=b_1$ on $E$. We shall also leave
$\exp_G$ implicit in the expression for $\Psi _{kj}$ and $A_{\varpi kj}$, which shall be viewed in the following
as functions of $\xi \in E$.
Thus, replacing $g$ by $\xi$ ad $dg$ by $H_G(\xi)\,d\xi$ in (\ref{eqn:j-th-summand-fourier-itegral}),
and then performing the change of variable $\xi=\nu/\sqrt{k}$, we obtain

\begin{eqnarray}\label{eqn:j-th-summand-fourier-itegral-exponential-chart}
\lefteqn{\Pi_{\varpi,k}\left
(x+\frac{w}{\sqrt{k}},x+\frac{v}{\sqrt{k}}\right)_{j}}\\
&\sim&k^{1-\mathrm{g}/2}\,e^{-ik\vartheta_j}\int_0^{+\infty} \,
\int _{-\epsilon}^{\epsilon}\,\int_{\mathbb{R}^\mathrm{g}}e^{ik \Psi _{kj}
\left(\frac{\nu}{\sqrt{k}},t,\vartheta\right)}\,A_{\varpi kj}\left(\frac{\nu}{\sqrt{k}},t,\vartheta\right)
\,H_G\left (\frac{\nu}{\sqrt{k}}\right )\,dt\,d\vartheta\,d\nu.\nonumber
\end{eqnarray}

Our next step will be to Taylor expand $\Psi _{kj}$ in descending powers of $k^{1/2}$, by relying
on (64) and (65) of \cite{sz}; to this end, we shall need the Heisenberg local coordinates of
$\mu _{g^{-1}}\circ r_{e^{i(\vartheta+\vartheta_j)}}\left( \mu_{g_j^{-1}}
\left (x+w/\sqrt{k}\right )\right)$ when $g=\exp_G\big(\nu/\sqrt{k}\big)$.

Recalling that $m=\pi(x)$ and $G_m\subseteq G$ is the stabilizer subgroup,
let us consider the isotropy representation
$G_m\rightarrow \mathrm{GL}\big(T_{m}M\big)$, $g\mapsto d_m\mu_g$;
for every $j=1,\ldots,N_x$ and $w\in T_mM$, let us set $w_j=:d_m\mu_{g_j^{-1}}(w)\in T_mM$.

In view of our choice of $\omega=\frac i2\,\Theta$ as the reference K\"{a}hler form
in our construction of Heisenberg local coordinates, we then have:

\begin{lem}
\label{lem:local coordinates} Suppose $x\in X$, $\Phi\circ \pi(x)=0$, and fix a system of
Heisenberg local coordinates centered at $x$. Then there exist
$\mathcal{C}^\infty$ functions $Q,T:\mathbb{C}^\mathrm{n}\times
\mathbb{R}^\mathrm{g}\rightarrow \mathbb{C}^\mathrm{n}$, vanishing
at the origin to third and second order, respectively, such that the
following holds. For every $w\in T_{\pi(x)}M$, $-\pi<\vartheta<\pi$,
$\nu \in \mathfrak{g}$, as $k\rightarrow +\infty$ the Heisenberg
local coordinates of
$$X_{j,k}(x,w,\nu)=:\mu_{e^{-\nu/\sqrt{k}}}\left (r_{e^{i\vartheta_j}}\circ
\mu_{g_j^{-1}}
\left (x+\frac{w}{\sqrt{k}}\right)\right )$$ are given by
\begin{eqnarray*}
\left (
\frac{1}{\sqrt{k}}\,\big [w_j-\nu_M(m)\big]+T_j\left
(\frac{w}{\sqrt{k}},\frac{\nu}{\sqrt{k}}\right),\frac 1k \,\omega _m\Big (\nu _M(m),w_j\Big)+Q_j\left
(\frac{w}{\sqrt{k}},\frac{\nu}{\sqrt{k}}\right)\right),\end{eqnarray*} where
$Q_j,T_j:\mathbb{C}^\mathrm{n}\times \mathbb{R}^\mathrm{g}\rightarrow
\mathbb{C}^\mathrm{n}$ vanish at the origin to third and second
order, respectively.
\end{lem}

\textit{Proof.} Set $m=\pi(x)$. By definition of $\nu_M$,
$\mu_{e^{-\nu/\sqrt{k}}}\circ \mu_{g_j^{-1}}\left (m+\frac{w}{\sqrt{k}}\right)$ has
underlying preferred local coordinates
$\frac{1}{\sqrt{k}}\,\big(w_j-\nu _M(m)\big)+T
\left(\frac{1}{\sqrt{k}}\,w_j,\frac{1}{\sqrt{k}}\,\nu \right)$,
where $T:\mathbb{C}^\mathrm{n}\times
\mathbb{R}^\mathrm{g}\rightarrow \mathbb{C}$ vanishes to second
order at the origin (here, $w,\,\nu_M(m)\in T_mM$ are identified
with their images in $\mathbb{C}^\mathrm{n}$ under the unitary
isomorphism $T_mM\rightarrow \mathbb{C}^\mathrm{n}$ induced by the
chosen preferred local coordinates centered at $m$).

Therefore, the Heisenberg local coordinates of
$X_{j,k}(x,w)$ have the form
$
\Big (\theta\left (k^{-1/2}\right),
k^{-1/2}\,\big(w_j-\nu_M(m)\big)+T\left
(k^{-1/2}\,w_j,k^{-1/2}\,\nu\right)\Big),$ for an
appropriate smooth function $\theta :(-\delta,\delta)\rightarrow
\mathbb{R}$.

In order to determine $\theta$, set $
\gamma _s(t)=:\mu_{e^{-t\nu}}\circ
\mu_{g_j^{-1}}\left (x+sw\right)$, defined and smooth for
all sufficiently small $s,t\in \mathbb{R}$.
Let us write
$w_j=p_{w_j}+iq_{w_j}$, $\nu _M(m)=p_\nu+iq_\nu$, where
$p_{w_j},q_{w_j},p_\nu,q_\nu\in \mathbb{R}^\mathrm{n}$. The preferred local
coordinates of $\pi \big(\gamma _s(t)\big)=\mu_{e^{-t\nu}}\circ
\mu_{g_j^{-1}}\left
(m+sw\right)$ are $(sp_{w_j}-tp_\nu)+i(sq_w-tq_\nu)+T(sw,t\nu).$
Therefore, the Heisenberg local coordinates of $\gamma _s(t)$ have
the form
\begin{equation}\label{eqn:local-form-for-gamma-s-t}
\Big(\widetilde{\theta}
(s,t),\big(sp_{w_j}-tp_\nu\big)+i\big(sq_{w_j}-tq_\nu\big)+T(sw,t\nu)\Big),
\end{equation} for an
appropriate smooth function $\widetilde{\theta}(s,t)$; clearly,
$\theta(u)=\widetilde{\theta}(u,u)$.

\begin{claim}\label{claim:form-of-theta}
We have
$
\widetilde{\theta}(s,t)=-\vartheta_j+(st)\cdot
d_0(w_j,\nu)+\widetilde{\theta}_1(sw_{j},t\nu)$, for appropriate  smooth
functions $d_0,\widetilde{\theta}_1:\mathbb{C}^\mathrm{n}\times
\mathbb{R}^\mathrm{g}\rightarrow \mathbb{R}$, with
$\widetilde{\theta}_1$ vanishing to third order at the origin.
\end{claim}

\textit{Proof of Claim \ref{claim:form-of-theta}.}
Recall that
Heisenberg local coordinates depend on the choice of a preferred
local holomorphic
frame $e_L$ of $L$
an open neighborhood $U\subseteq M$ of $m$; as discussed in \S \ref{sctn:prelimiaries},
without loss of generality we may assume
that $U$ is $G_m$-invariant and $g^*(e_L^*)=h_g\cdot e_L^*$, $\forall\,g\in G_m$.
Let us write $\sigma=e_L^*$.
We have $x+sw=\sigma (m+sw)/\|\sigma (m+sw)\|$, where $m+sw\in U$ is the point with local preferred
holomorphic coordinates $w\in \mathbb{C}^\mathrm{n}$.
Therefore,
\begin{eqnarray*}
\lefteqn{\mu_{g_j}^{-1}(x+sw)=\mu_{g_j}^{-1}\Big(\sigma (m+sw)\Big)\Big/\|\sigma (m+sw)\|}\\
&=&h _{g_j}^{-1}\,
\sigma\Big(\mu_{g_j}^{-1}(m+sw)\Big)\Big/\left\|\sigma\Big(\mu_{g_j}^{-1}(m+sw)\Big)\right\|
\end{eqnarray*}
has Heisenberg local coordinates
$\Big(-\vartheta_j,z_j(w,s)\Big)$, where $z_j(w,s)$ are the local preferred
holomorphic coordinates of $\mu_{g_j}^{-1}(m+sw)$. Therefore,
$\widetilde{\theta}
(s,0)=-\vartheta_j$ for all $s$.

We conclude that
$\widetilde{\theta}(s,t)=-\vartheta_j+t\,R(s,t)$, for some smooth function $R$.

On the other hand, $X'$ is $G$-invariant, and $G$ acts horizontally
on it (in other words, for every $x\in X'$ and $\xi \in
\mathfrak{g}$ we have $\xi _X(x)=\xi _M^\sharp\big (\pi(x)\big)$,
where $\xi _M^\sharp$ denotes the horizontal lifting of $\xi _M$).
Lemmata 2.4 and 3.3 of \cite{dp} then imply that $\theta (0,t)=t^3\,S(t)$ for a smooth function $S(t)$.
Thus, $R(s,t)=t^2 R_1(t)+s\,d(s,t)$ for smooth functions $R_1(t)$, $d(s,t)$.
We conclude that $\widetilde{\theta}(s,t)=t^3\,R_1(t)+st\,d(s,t)$, and the statement follows by setting
$d_0=d(0,0)$.

\hfill Q.E.D.

\medskip

Returning to the proof of Lemma \ref{lem:local coordinates},
in order to determine $d_0$ we recall that the expression for $\alpha$
in Heisenberg local coordinate
is $\alpha=d\theta+p\,dq-q\,dp+\beta(\|z\|^2)$.
Inserting the local expression for $\gamma _s(t)$ that we obtain from (\ref{eqn:local-form-for-gamma-s-t})
and Claim \ref{claim:form-of-theta},
we obtain
\begin{eqnarray}\label{eqn:espressione-locale-per-alpha}
\lefteqn{\gamma _s^*(\alpha)}\nonumber\\
&=&\left\{\Big [sd_0+
\big(sp_{w_j}-tp_\nu\big)\cdot (-q_\nu)-\big(sq_{w_j}-tq_\nu\big)\cdot (-p_\nu)dt\Big]+G_1(s,t,\nu,w_{j})\right \}\,\,dt \nonumber\\
&=&\left \{s\,\Big [d_0+\omega_m(\nu_M,w_{j})\Big ]+G_1(s,t,\nu,w_{j})\right \}\,dt,
\end{eqnarray}
where $G_1$ vanishes to second order for $(s,t)=(0,0)$.

On the other hand, we have $d_t(\gamma _s)(1)=-\nu _X\big(\gamma _s(t)\big)$. Therefore
\begin{equation}\label{eqn:gamma-con-alfa}
\gamma _s^*(\alpha)(t)=-\iota(\nu_X)\alpha \Big(\pi\big(\gamma _s(t)\big)\Big)\,dt
=\Phi ^\nu\Big (\pi\big(\gamma _s(t)\big)\Big)\,dt,
\end{equation}
having used that $\Phi^\nu=-\iota (\nu_X)\alpha$.

Because of the $G$-equivariance of $\Phi$,
$\Phi\circ \pi\Big(\gamma _0(t)\Big)=\Phi\circ \pi\Big(\mu_{g_je^{t\nu}}(x)\Big)=0$ for every sufficiently small $t$; therefore,
$\left .\frac{\partial \Phi \circ \gamma}{\partial t}\right |_{(0,t)}=0$ identically, where with
abuse of language we have written $\Phi$ for $\Phi\circ \pi:X\rightarrow \mathfrak{g}$.
This implies
\begin{eqnarray}\label{eqn:phi-vicino}
\Phi ^\nu\Big (\pi\big(\gamma _s(t)\big)\Big)&=&s\,d_m\Phi^\nu(w_j)+G_3(s,t,\nu,w_j)\nonumber\\
&=&2s\,\omega_m(\nu_M,w_j)+G_3(s,t,\nu,w_j),
\end{eqnarray}
where $m=\pi(x)$, and $G_3$ vanishes to second order for $(s,t)=(0,0)$.

Comparing (\ref{eqn:gamma-con-alfa}) and (\ref{eqn:phi-vicino})
with (\ref{eqn:espressione-locale-per-alpha}), we obtain $d_0=\omega_m(\nu_M,w)$, $G_2=G_3$.
To complete the proof of Lemma \ref{lem:local coordinates}, we need only take $s=t=1/\sqrt{k}$, and remark that
in Heisenberg local coordinates $r_{e^{i\vartheta_j}}$ is simply translation by $\vartheta_j$.

\hfill Q.E.D.

\medskip

Let us set $\psi_2(u,v)=u\cdot \overline{v}-\frac 12\left (\|u\|^2+\|v\|^2\right)$
($u,\,v\in \mathbb{C}^\mathrm{n}$).
Invoking (63)-(65) of \cite{sz}, in view of Lemma \ref{lem:local coordinates}
we obtain that $\Psi_{kj}$ in
(\ref{eqn:psi-j-k}) has the form:

\begin{eqnarray}\label{eqn:psi-kj}
\lefteqn{\Psi_{kj}(g,t,\vartheta)}\\
&=&i\,t\,\left (1-e^{i\vartheta}\right)-\vartheta+\frac{t}{k}\,e^{i\vartheta}\,
\Big[\omega _m\big(\nu_M(m),w_j\big)-i\,\psi_2\big(w_j-\nu_M(m),v\big)\Big]\nonumber \\
&&+t\,e^{i\vartheta}\,R_j\left (\frac{\nu_M(m)}{\sqrt{k}},\frac{w}{\sqrt{k}},\frac{v}{\sqrt{k}}\right),
\nonumber\end{eqnarray}
where $R_j:\left (\mathbb{C}^\mathrm{\mathrm{n}}\right )^3\rightarrow \mathbb{C}$ is a smooth function
vanishing to third order at
the origin.

Let us now insert (\ref{eqn:psi-kj}) in (\ref{eqn:j-th-summand-fourier-itegral-exponential-chart}).
We obtain

\begin{eqnarray}\label{eqn:j-th-summand-fourier-itegral-exponential-chart-fixed-phase}
\lefteqn{\Pi_{\varpi,k}\left
(x+\frac{w}{\sqrt{k}},x+\frac{v}{\sqrt{k}}\right)_{j}}\\
&\sim&k^{1-\mathrm{g}/2}\,e^{-ik\vartheta_j}\int_0^{+\infty} \,
\int _{-\epsilon}^{\epsilon}\,\int_{\mathbb{R}^\mathrm{g}}e^{ik \Psi
\left(t,\vartheta\right)}\,\widetilde{A}_{\varpi kj}\left(\nu,w,v,t,\vartheta\right)
\,dt\,d\vartheta\,d\nu,\nonumber
\end{eqnarray}
where we have set
\begin{equation}\label{defn:fixed-phase}
\Psi
\left(t,\vartheta\right)=:i\,t\,\left (1-e^{i\vartheta}\right)-\vartheta,\end{equation}
\begin{eqnarray}\label{defn:final-amplitude}
\lefteqn{
\widetilde{A}_{\varpi kj}
\left(\nu,w,v,t,\vartheta\right)=:e^{t\,e^{i\vartheta}\,\big[\psi_2\big(w_j-\nu_M(m),v\big)+
i\omega _m\big(\nu_M(m),w_j\big)\big]}
}\\
&&\cdot e^{i\,k\,t\,e^{i\vartheta}\,R_j\left (\frac{\nu_M(m)}{\sqrt{k}},\frac{w}{\sqrt{k}},\frac{v}{\sqrt{k}}\right)}\,\,A_{\varpi kj}\left(\frac{\nu}{\sqrt{k}},t,\vartheta\right)
\,H_G\left (\frac{\nu}{\sqrt{k}}\right ).\nonumber
\end{eqnarray}

A straightforward computation shows that
$$
\psi_2\big(w_j-\nu_M(m),v\big)+
i\omega _m\big(\nu_M(m),w_j\big)=T_\mathrm{h}+T_\mathrm{t}+T_\mathrm{v}+T_\mathrm{vt},
$$
where
$$T_\mathrm{h}=:\psi_2\big(w_{\mathrm{h}},v_\mathrm{h}\big),\,\,\, T_\mathrm{t}=:-\frac 12\,\big\|w_{\mathrm{t}}-
v_\mathrm{t}\big\|^2, \,\,\,T_\mathrm{v}=:-\frac 12\,\big\|w_{j\mathrm{v}}-\nu_M(m)-v_{j\mathrm{v}}\big\|^2,$$
$$T_\mathrm{tv}=:-i\,\omega _m\big(w_{j\mathrm{v}}-\nu_M(m)-v_{j\mathrm{v}},w_{j\mathrm{t}}+v_{j\mathrm{t}}\big)+
i\Big[\omega _m\big(w_\mathrm{v},w_\mathrm{t}\big)-\omega _m\big(v_\mathrm{v},v_\mathrm{t}\big)\Big]$$
(notice that the map $w\mapsto w_j$ induced by the isotropy action of
$g_j^{-1}\in G_m\subseteq G$ is an isometry of $T_mM$, since $G$ preserves the metric of $M$).

We may insert in (\ref{eqn:a-j-k})
the asymptotic expansion for the classical symbol $s(x,y,t)$ appearing in the parametrix
for $\Pi$, and use Taylor expansion in $g=\nu/\sqrt{k}$, $w/\sqrt{k}$ and
$v/\sqrt{k}$ in descending powers of
$k^{1/2}$, to deduce that

\begin{eqnarray}\label{defn:final-amplitude-asymptotic-expansion}
\widetilde{A}_{\varpi kj}
\left(\nu,w,v,t,\vartheta\right)\sim \sum _{l\ge 0}a_{\varpi jl}
\left(\nu,w,v,t,\vartheta\right)\,k^{\mathrm{n}-l/2},
\end{eqnarray}
where every coefficient has the form
$$a_{\varpi jl}
\left(\nu,w,v,t,\vartheta\right)=e^{te^{i\vartheta}(T_\mathrm{h}+T_\mathrm{t}+T_\mathrm{v}+T_\mathrm{vt})}\,
p_{\varpi jl}
\left(\nu,w,v,t,\vartheta\right)$$ and each $p_{\varpi jl}
\left(\nu,w,v,t,\vartheta\right)$
is a polynomial in $\nu$, $w$ and $v$ with coefficients depending on $x,t,\vartheta$ and $\varpi$.
In particular, the leading coefficient is

\begin{eqnarray}\label{defn:final-amplitude-asymptotic-expansion-leading-term}
a_{\varpi j0}
\left(\nu,w,v,t,\vartheta\right)= e^{
t\,e^{i\vartheta}
\,\Big(T_\mathrm{h}+T_\mathrm{t}+T_\mathrm{v}+T_\mathrm{vt}\Big)}
\,\gamma _0(\vartheta)\,\beta_\varpi(g_j)\,s_0(x,x)\,t^\mathrm{n}.
\end{eqnarray}

Thus,
\begin{eqnarray}\label{eqn:j-th-summand-fourier-itegral-exponential-expanded}
\lefteqn{\Pi_{\varpi,k}\left
(x+\frac{w}{\sqrt{k}},x+\frac{v}{\sqrt{k}}\right)_{j}}\\
&\sim&k^{1-\mathrm{g}/2}\,e^{-ik\vartheta_j}\,\sum_{l\ge 0}\left(\int_0^{+\infty} \,
\int _{-\epsilon}^{\epsilon}\,\int_{\mathbb{R}^\mathrm{g}}e^{ik \Psi
\left(t,\vartheta\right)}\,a_{\varpi jl}
\left(\nu,w,v,t,\vartheta\right)\,k^{\mathrm{n}-l/2}
\,dt\,d\vartheta\,d\nu\right),\nonumber
\end{eqnarray}

To determine the leading asymptotics of (\ref{eqn:j-th-summand-fourier-itegral-exponential-expanded}),
let us first integrate (\ref{defn:final-amplitude-asymptotic-expansion-leading-term})
in $d\nu$.
By our choice of Heisenberg local coordinates, we may unitarily identify
$(T_mM,\omega_m)$ with $(\mathbb{C}^\mathrm{n},\omega_0)$, where
$\omega_0$ is the standard symplectic structure on $\mathbb{C}^\mathrm{n}$; let $g_0$ be the
standard scalar product on $\mathbb{C}^\mathrm{n}$, so that
$\omega _0\big(\mathbf{a},\mathbf{b}\big)=-g_0\big(\mathbf{a},J_0(\mathbf{b})\big)$,
$\forall\,\mathbf{a},\,\mathbf{b}\in \mathbb{C}^\mathrm{n}$, where $J_0$ is multiplication by $i$.
We shall view $S_m:\nu\mapsto \nu _M(m)$
as a map $\mathfrak{g}\rightarrow \mathbb{C}^\mathrm{n}$.
Let us set $\lambda =
t\,e^{i\vartheta}$. Up to a multiplicative factor, we are led to integrating
\begin{equation}\label{eqn:integral-for-leading-coeff}
e^{\lambda \,\left [-\frac 12\,\big\|w_{j\mathrm{v}}-\nu_M(m)-v_{j\mathrm{v}}\big\|^2+i\,
g_0\big(w_{j\mathrm{v}}-\nu_M(m)-v_{j\mathrm{v}},J_0(w_{j\mathrm{t}}+v_{j\mathrm{t}})\big)\right]}
\end{equation}
in $d\nu$.

Recall that the $\nu$ coordinates are induced by the choice of an orthonormal
basis of $\mathfrak{g}$; we can shift the integration to the tangent space of the $G$-obit
through $m$, $\mathfrak{g}_M(m)\subseteq T_mM$. Let us then choose an orthonormal basis of
$\mathfrak{g}_M(m)$, and let $\beta$ be the corresponding linear coordinates.
We can use $\beta$ as integration variable, by the relation
$\beta=S_m(\nu)$. By Lemma 3.9 of \cite{dp},
after performing the change of variables $\beta\mapsto \beta-(w_{j\mathrm{v}}-v_{j\mathrm{v}})$ we are
left with
\begin{eqnarray}\label{eqn:valutazione-integrale-dnu}
\lefteqn{\frac{1}{V_{\mathrm{eff}}(m)\,|G_m|}\,\int _{\mathbb{R}^\mathrm{g}}\,
e^{t\,e^{i\vartheta}\,\left[-\frac 12\|\beta\|^2-
i\,g_0\big(\beta,J_0(w_{j\mathrm{t}}+v_{j\mathrm{t}})\big)\right]}\,d\beta}\\
&=&\frac{(2\pi)^{\mathrm{g}/2}}{V_{\mathrm{eff}}(m)\,|G_m|}\cdot
\frac{1}{\sqrt{t}\,e^{i\vartheta/2}}\,e^{-\frac 12\,t\,e^{i\vartheta}\,\|w_{\mathrm{t}j}+v_{\mathrm{t}j}\|^2};
\nonumber
\end{eqnarray}
in fact, since $t>0$
(\ref{eqn:valutazione-integrale-dnu}) is valid when $\vartheta=0$ because
$-\frac 12\|\beta\|^2$ equals its own Fourier transform, and consequently by analytic
continuation it holds for all $\vartheta \in (-\pi/2,\pi/2)$.

Let us next consider the case a general $a_{\varpi jl}
\left(\nu,w,v,t,\vartheta\right)$. Up to multiplicative factors polynomial in $w$ and $v$,
we are led to integrate the product of (\ref{eqn:integral-for-leading-coeff})
times a monomial in $\nu$.
Again up to an appropriate scalar factor, this amounts to multiplying the integrand in (\ref{eqn:valutazione-integrale-dnu})
by a monomial in $\beta$, hence to evaluating an appropriate higher
derivative of $e^{-\|\beta\|^2/2}$ in
$J_0(w_{j\mathrm{t}}+v_{j\mathrm{t}})\in \mathfrak{g}_M(m)$.
We are thus left with the product of the right hand side in
(\ref{eqn:valutazione-integrale-dnu})
times a polynomial in $w_t$ and $v_t$.

We can now insert (\ref{eqn:valutazione-integrale-dnu}) in (\ref{defn:final-amplitude-asymptotic-expansion})
and (\ref{eqn:j-th-summand-fourier-itegral-exponential-chart-fixed-phase}) to obtain

\begin{eqnarray}\label{eqn:j-th-summand-fourier-itegral-exponential-in-t-theta}
\lefteqn{\Pi_{\varpi,k}\left
(x+\frac{w}{\sqrt{k}},x+\frac{v}{\sqrt{k}}\right)_{j}}\\
&\sim&k^{1-\mathrm{g}/2}\,e^{-ik\vartheta_j}\int_0^{+\infty} \,
\int _{-\epsilon}^{\epsilon}\,e^{ik \Psi
\left(t,\vartheta\right)}\,S_{\varpi j}\left(w,v,t,\vartheta,k\right)
\,dt\,d\vartheta\,\nonumber
\end{eqnarray}
where $S_{\varpi j}\left(w,v,t,\vartheta,k\right)\sim \sum_{l\ge 0} S_{\varpi j l}(w,v,t,\vartheta)\,
k^{\mathrm{n}-l/2}$ and the coefficients of the expansion are as follows.

First,
the leading coefficient is

\begin{equation}
S_{\varpi j 0}(w,v,t,\vartheta)=\frac{(2\pi)^{\mathrm{g}/2}}{V_{\mathrm{eff}}(m)\,|G_m|}\cdot
\frac{1}{\sqrt{t}\,e^{i\vartheta/2}}\,\gamma _0(\vartheta)\,\beta_\varpi(g_j)s_0(x,x)\,t^\mathrm{n}\,
e^{t\,e^{i\vartheta}\,\Gamma (w,v)}
\end{equation}
where $\Gamma (w,v)=\psi_2(w_\mathrm{h},v_\mathrm{h})-\big\|w_{\mathrm{t}}\big\|^2
-\big\|v_{\mathrm{t}}\big\|^2+i\Big[\omega _m\big(w_\mathrm{v},w_\mathrm{t}\big)-\omega _m\big(v_\mathrm{v},v_\mathrm{t}\big)\Big]$.

Next, for every $l\ge 1$ we have $S_{\varpi j l}(w,v,t,\vartheta)=p_{\varpi j l}(w,v,t,\vartheta)\,
e^{t\,e^{i\vartheta}\,\Gamma (w,v)}$, where $p_{\varpi j l}(w,v,t,\vartheta)$ is a poynomial
in $w$ and $v$.

Thus
we are left with an oscillatory integral whose phase $\Psi$, given by
(\ref{defn:fixed-phase}), is the same phase appearing in
the discussion of the scaling asymptotics of non-equivariant Szeg\"{o} kernels in \S 3 of \cite{sz}.
In particular, $\Psi$ has non-negative imaginary part, and a unique stationary point for $t=1$ and $\vartheta=0$;
furthermore, at this point the Hessian of $\Psi$ is
$$
\Psi''(1,0)=\left(
              \begin{array}{cc}
                0 & 1 \\
                1 & i \\
              \end{array}
            \right).
            $$
Hence $(1,0)$ is a non-degenerate stationary point of $\Psi$.
Arguing as in \textit{loc. cit.}, the contribution coming from $|t|\ge 2$, say,
is rapidly decreasing, and by the stationary phase method
for complex oscillatory integrals (Theorem 7.7.5 of \cite{hor})
there is an asymptotic expansion:
\begin{eqnarray}\label{eqn:j-th-summand-fourier-itegral-exponential-in-t-theta-stationary-phase}
\lefteqn{\Pi_{\varpi,k}\left
(x+\frac{w}{\sqrt{k}},x+\frac{v}{\sqrt{k}}\right)_{j}}\\
&\sim&k^{1-\mathrm{g}/2}\,e^{-ik\vartheta_j}\frac{1}{\sqrt{\det\big(k\Psi''(1,0)/2\pi i\big)}}\,\sum_{s=0}
^{+\infty}\,k^{-s}\,\left.L_s\big(S_{\varpi j}\left(w,v,t,\vartheta,k\right)\big)\right|_{t=1,\vartheta=0},\nonumber
\end{eqnarray}
where $L_0$ is the identity, and $L_s$ is a suitable differential operator of degree $2s$ in $(t,\theta)$
for any $s=0,1,2\ldots$.
The statement then follows from the previous description of the phase; in particular, each coefficient
in the asymptotic expansion is the product of $e^{\Gamma (w,v)}$
and a polynomial in $w,\,v$.

The statement of the Theorem follows by summing over $j$.

\end{document}